\documentstyle[12pt]{article}

\newcommand{\sect}[1]{\setcounter{equation}{0}\section{#1}}
\newcommand{\subsect}[1]{\subsection{#1}}

\newtheorem{proposition}{Proposition}

\def\be{\begin{equation}}
\def\ee{\end{equation}}
\def\bea{\begin{eqnarray}}
\def\eea{\end{eqnarray}}

\def\te{\chi}

\def\m{\mu}
\def\Ck{{\rm{\ \!C\/}}}
\def\Sk{{\rm{\ \!S\/}}}

\def\>#1{{\bf #1}}
\def\1{\'{\i}}
\def\R{{\rm I\kern-.2em R}}
\def\back{\!\!\!\!\!\!}

\def\tfrac#1#2{{\scriptstyle{\frac{#1}{#2}}}}

\begin{document}

\thispagestyle{empty}

\hfill\today

\vspace{2cm}

\begin{center}
{\LARGE{\bf{Universal $R$--matrices for }}}

{\LARGE{\bf{ non-standard (1+1) quantum  groups}}}

\end{center}

\bigskip\bigskip\bigskip

\begin{center}
A. Ballesteros$^{1}$, E. Celeghini$^{2}$, F.J. Herranz$^{1}$, M.A. del
Olmo$^{3}$ and M. Santander$^{3}$
\end{center}

\begin{center}
\em

{ {  $^{1}$ Departamento de F\1sica, Universidad de Burgos}
\\  E--09003, Burgos, Spain}

{ $^{2}$ Dipartimento di Fisica and INFN, Sezione di
Firenze\\ Largo E. Fermi, 2. I--50125 Firenze, Italy}

{ { $^{3}$ Departamento de F\1sica Te\'orica, Universidad de Valladolid} \\
E--47011, Valladolid, Spain }

   \end{center}
\rm

\bigskip\bigskip\bigskip

\begin{abstract}
A universal quasitriangular $R$--matrix for the non-standard  quantum
(1+1) Poincar\'e algebra $U_ziso(1,1)$ is deduced  by imposing
analyticity in the deformation parameter $z$. A family $g_\mu$ of ``quantum
graded contractions" of the algebra $U_ziso(1,1)\oplus U_{-z}iso(1,1)$ is
obtained; this set of quantum algebras contains as Hopf subalgebras with two
primitive translations quantum analogues of the  two dimensional Euclidean,
Poincar\'e and Galilei algebras enlarged  with dilations. Universal
$R$--matrices
 for these quantum Weyl algebras and their associated quantum groups are
constructed.

\end{abstract}

\newpage

\sect {Introduction}

Two types of  quantum deformations for the $so(2,2)$ algebra
and for its most relevant graded contractions have been  recently
studied in   \cite{Beyond}. They are called {\it standard} and
{\it non-standard} quantum algebras according to the fact that their
corresponding  coboundary Lie bialgebras  come from  a classical
$r$--matrix which is a skew solution either of the modified  classical
Yang--Baxter equation (YBE) or of the  classical YBE, respectively.
In contradistinction with
the standard case, the family of non-standard  quantum  algebras
contains as quantum Hopf subalgebras
  two dimensional Euclidean and (1+1)
Poincar\'e and Galilei algebras enlarged with a dilation:  the
so-called  ``Weyl" or similarity subalgebras. These quantum subalgebras share
the property of including the two translation generators as
primitive. This fact could be  relevant  in relation with the problem
of discretizing two dimensional spaces in some symmetric way.

Let us also recall that a quasitriangular Hopf algebra \cite{Dr} is  a
pair $({{\cal A}},{{\cal R}})$
 where ${\cal A}$ is a Hopf algebra and ${\cal R}  \in {\cal A}\otimes
{\cal A}$ is invertible and verifies
\bea
&&\sigma \circ \Delta  X = {\cal R} (\Delta  X) {\cal R}^{-1},
 \qquad \forall\
\, X\in {\cal A}, \label{aa}\\
&& (\Delta  \otimes id){\cal R} ={\cal R}_{13}{\cal
R}_{23},\qquad   (id \otimes \Delta ){\cal R} ={\cal R}_{13}{\cal R}_{12},
\label{ab}
\eea
where, if ${\cal R}=\sum_{i}
a_i\otimes b_i$,  we denote ${\cal R}_{12}\equiv\sum_{i} a_i\otimes
b_i\otimes 1$, ${\cal R}_{13}\equiv\sum_{i} a_i\otimes 1\otimes b_i$,
${\cal R}_{23}\equiv\sum_{i} 1\otimes a_i\otimes b_i$ and $\sigma$ is
the flip operator $\sigma(x\otimes y)=(y\otimes x)$. If ${\cal A}$ is a
quasitriangular Hopf algebra, then  ${\cal R}$ is called a universal
${\cal R}$--matrix and satisfies the quantum  YBE:
\be
{\cal R}_{12}{\cal R}_{13}{\cal R}_{23}
={\cal R}_{23}{\cal R}_{13}{\cal R}_{12}.\label{ac}
\ee

In this paper,
the construction of universal $R$--matrices and quantum groups for
the above mentioned
family of non-standard  quantum algebras is discussed.   In particular,
we obtain the universal $R$--matrices for the  Weyl subalgebras.

A  straightforward approach to this problem would consist on
starting from a universal $R$--matrix for the non-standard Hopf
algebra $U_zsl(2,\R)\simeq U_zso(2,1)$, since the prescription
$U_zso(2,2) \simeq U_zso(2,1)\oplus
U_{-z}so(2,1)$ applied to $R$--matrices would lead to a universal
$R$--matrix for $U_zso(2,2)$. Afterwards, by introducing   ``quantum
graded contractions" a set of $R$--matrices for all the family of
non-standard algebras would be obtained. A similar procedure was
developed in \cite{CGST1,CGST2} for getting universal $R$--matrices for
some standard quantum algebras.  Unfortunately, to our
knowledge, no universal  $R$--matrix for the non-standard
$U_zsl(2,\R)$ has appeared
in the literature, since in fact the one given in \cite{Ohn}  neither
verifies (\ref{aa}) nor satisfies (\ref{ac}).

Therefore, we have to focus the problem from a different point of view.
Following the method developed in  \cite{Heis2} for the standard (1+1)
groups and in \cite{Heis} for the Heisenberg group,
  we impose   analyticity in
the deformation parameter $z$ and relation (\ref{aa})
in order to get an $R$--matrix for the non-standard
quantum (1+1) Poincar\'e algebra $U_ziso(1,1)$. The $R$--matrix  so
obtained coincides in turn with a universal $R$--matrix for the
positive Borel subalgebra of the non-standard $sl(2,\R)$ given in
\cite{AA}; this fact proves its universality. This explicit
construction together with a brief overview of both the quantum Poincar\'e
algebra and group studied in \cite{Tmatrix}
  is  presented in Section 2.  The   fact that
\be
U_zt_4(so(1,1)\oplus so(1,1))\simeq U_ziso(1,1)\oplus
U_{-z}iso(1,1)\label{ad}
\ee
 can be also used at the group and   $R$--matrix
levels
 and leads to the whole quantum structure for this group as it is shown
in Section 3; Lie algebras with structure $t_n(so(p,q)\oplus so(r,s))$
are described in \cite{BW}. A ``quantum graded contraction" is introduced
providing quantum structures for   two more non-standard quantum real
algebras: $U_ziiso(1,1)$ and  $U_zt_4(so(2)\oplus so(1,1))$ (the latter is
isomorphic to the (2+1) expanding Newton--Hooke algebra). Each   of
these quantum algebras  contains a    Hopf Weyl subalgebra. Section 4
is devoted to the obtention of
the   universal $R$--matrices and quantum groups corresponding to these
Weyl  subalgebras.

\sect {Universal $R$--matrix for the  Poincar\'e group}

 The (1+1) Poincar\'e algebra $iso(1,1)$  is generated by
one boost generator $K$ and the translation generators
along the light-cone $P_\pm$.
The Lie brackets are:
\be
[K,P_\pm]=\pm 2 P_\pm,\qquad [P_+,P_-]=0 .\label{ba}
\ee
 A non-standard coboundary bialgebra of $iso(1,1)$ is generated by
the classical $r$--matrix $r=z K\wedge P_+$ which verifies the
classical YBE.

The quantum deformations for the universal enveloping algebra
$Uiso(1,1)$ and for the algebra of smooth functions on the group
$Fun(ISO(1,1))$, denoted respectively by $U_ziso(1,1)$ and
$Fun_z(ISO(1,1))$, are given by the following statements (see
 \cite{Tmatrix} for a more detailed exposition and
proofs and also \cite{Lyak}):

\begin{proposition}
 The Hopf structure of  $U_ziso(1,1)$ is given by the  coproduct,
counit, antipode
\bea
&&\Delta P_+ =1 \otimes P_+ + P_+\otimes 1,\nonumber\\
&&\Delta P_- =e^{-zP_+} \otimes P_- + P_-\otimes e^{zP_+},\label{bb}\\
&&\Delta K =e^{-zP_+} \otimes K + K\otimes e^{zP_+};\nonumber
\eea
\be
\epsilon(X) =0;\qquad
\gamma(X)=-e^{zP_+}\ X\ e^{-zP_+},\ \ \mbox{for $X\in \{K,P_\pm\}$};
\label{bc}
\ee
and the commutation relations
\be
[K,P_+]=2\,\frac{\sinh  zP_+ }z,\quad [K,P_-]=- 2\, P_-\cosh  zP_+  ,
\quad  [P_+,P_-]= 0 .\label{bd}
\ee
\end{proposition}

\begin{proposition}
 The  Hopf algebra
$Fun_z(ISO(1,1))$  has multiplication given by
\be
[{\hat\te},{\hat a_+}]=z(e^{ 2\hat \te} -1) ,\quad
[{\hat \te},{\hat a_-}]=0,\quad [{\hat a_+},{\hat a_-}]=- 2z \hat  a_-
;\label{bf} \ee
 coproduct
\be
\Delta(\hat \te)=\hat \te\otimes 1 + 1 \otimes\hat  \te,\qquad
\Delta(\hat a_\pm)=\hat a_\pm\otimes 1 + e^{\pm 2 \hat   \te}\otimes\hat
a_\pm  ;\label{bg} \nonumber
\ee
counit and antipode
\be
\epsilon(X) =0,\qquad
X\in\{\hat a_+,\hat a_-,\hat  \te\};
\label{bh}
\ee
\be
\gamma(\hat  \chi)=-\hat  \chi,\quad
 \gamma(\hat  a_\pm)=- e^{\mp 2\hat \te}\hat  a_\pm  .
\label{bi}
\ee
\end{proposition}

The quantum coordinates $\hat \chi$,  $\hat a_-$ and $\hat a_+$ of
$Fun_z(ISO(1,1))$  are, respectively, the  dual basis of the
 $U_ziso(1,1)$ generators $H=e^{zP_+}K$, $A_-=e^{-zP_+}P_-$ and
$A_+=P_+$.

Now we proceed to deduce a universal $R$--matrix for $U_ziso(1,1)$. We
assume that   the $R$--matrix is analytical in the quantum
parameter $z$ $(=\log q)$ and that   $R=1\otimes 1 +
z K\wedge P_+ + o(z^2)$. Hence we consider  an
$R$--matrix as a formal power series in $z$ with coefficients in $U
iso(1,1)\otimes U iso(1,1)$. We start from the following Ansatz:
\be
R=\exp \{z f(K,P_+,z)\, g(P_+,P_-,z)\}.\label{bj}
\ee
Firstly, we impose $R$ to verify relation  (\ref{aa}). Starting with the
primitive generator $P_+$, it is implied that
\be
[R,\Delta P_+]=0.\label{bk}
\ee
This requirement is fulfilled if
\be
[f(K,P_+,z)\, g(P_+,P_-,z), \Delta P_+]=
[f(K,P_+,z), \,\Delta P_+]\, g(P_+,P_-,z)=0.\label{bl}
\ee
Therefore, by taking into account commutation rules (\ref{bd}) a
solution for  $f(K,P_+,z)$ is
\be
f(K,P_+,z)=K\wedge \sinh zP_+.\label{bm}
\ee
Now we should   apply condition (\ref{aa}) for the two remaining
generators $P_-$ and $K$. Omitting the arguments of the functions $f$
and $g$ we have: \bea
R\, \Delta  X\, R^{-1}&\!=\!&\exp(zfg )
\Delta  X \exp(- zfg)\cr
&\!=\!&\Delta  X + z\,  [fg,\Delta X] + \frac {z^2}{2!}\,[fg,[fg,\Delta
X]] +\dots\cr &&\quad +\frac {z^n}{n!} \,
[fg,[fg\dots [fg,\Delta X]^{n)}\dots
]]+\dots \label{bn}
\eea
Since $P_-$ commutes with $P_+$, relation (\ref{bn}) with $X\equiv P_-$
becomes into
\bea
&&\exp(zfg )
\Delta  P_- \exp(- zfg)=\Delta P_- + z\,
[f,\Delta P_-]\,g +\frac {z^2}{2!}\,[f
,[f ,\Delta P_-]]\,g^2 +\dots\cr
&&\qquad\qquad \qquad\qquad\qquad +\frac {z^n}{n!}\,  [f ,[f\dots [f
,\Delta P_-]^{n)}\dots ]]\,g^n+\dots \label{bo}
\eea
We need to obtain the brackets $[f,\Delta P_-]$,
$[f,[f,\Delta P_-]]$,\dots\ in (\ref{bo}). The first and the second ones
are
\be
 [f,\Delta P_-]= A ,\qquad [f,[f,\Delta P_-]]=B\, 2\sinh z\Delta P_+  ,
\label{bp}
\ee
where
\bea
&& A\equiv 2 \exp(-z\Delta P_+)\sinh z P_+\otimes P_-
 -  2 \exp(z\Delta P_+) P_-\otimes \sinh zP_+,\\
&& B\equiv 2 \exp(-z\Delta P_+)\sinh z P_+\otimes P_-
 +  2 \exp(z\Delta P_+) P_-\otimes \sinh zP_+ .
\eea
Due to expression (\ref{bm}) $f$ commutes with any arbitrary function
 of $P_+$, so we
get
\be
[f,\sinh z\Delta P_+]=0,\qquad [f,B]= A\, 2\sinh z\Delta P_+ .
\ee
By a recurrence method we obtain the $2n$ and $2n+1$
iterates:
\bea
&&[f,[f \dots [f,\Delta P_-]^{2n)}\dots ]] = B\,  (2\sinh z\Delta
P_+)^{2n-1},\\ &&[f,[f \dots [f,\Delta P_-]^{2n+1)}\dots ]] = A\,
(2\sinh z\Delta P_+)^{2n} , \eea
so that expression (\ref{bo}) can be written as follows
\bea
&&\back\back\back\exp(zfg )  \Delta  P_- \exp(- zfg )  \cr
&&\back\back\back
= \Delta P_- + A\sum_{l=0}^\infty \frac {z^{2l+1}}{(2l+1)!} (2\sinh
z\Delta P_+)^{2l} g^{2l+1}  +
B\sum_{l=1}^\infty \frac {z^{2l}}{(2l)!}(2\sinh z\Delta P_+)^{2l-1}
g^{2l}\cr
&&\back\back\back = \Delta  P_- + \frac {A}{2\sinh z\Delta P_+}
\sinh(2z\sinh( z\Delta P_+)g)\cr
&&
+\frac {B}{2\sinh z\Delta P_+}\left[\cosh(2z\sinh
 (z\Delta P_+)g)-1\otimes 1\right].
\label{bba}
\eea
By introducing
\bea
&& C=\frac {1}{2\sinh z\Delta P_+}\sinh(2z\sinh(
z\Delta P_+)g),\\
&&D=\frac {1}{2\sinh z\Delta P_+}\left[\cosh(2z\sinh (z\Delta
P_+)g)-1\otimes 1\right] , \eea
expression (\ref{bba}) can be written as $\Delta  P_- + A C + B D$,
which
 must be equal to
\be
\sigma\circ \Delta P_- =e^{zP_+}\otimes P_- + P_- \otimes e^{-zP_+}.
\label{bbb} \ee
Thus, if we impose that (\ref{bba})  coincides with (\ref{bbb})    we
obtain the following system of equations for the function $g$:
\bea
&&\back\back\back
(1\otimes P_-)(e^{-zP_+}\otimes 1 - e^{zP_+}\otimes 1 + 2e^{-zP_+}
\sinh zP_+\otimes
e^{-zP_+} (D+C))=0,\\
&&\back\back\back
(P_-\otimes 1)(1\otimes e^{zP_+} - 1\otimes e^{-zP_+} + 2e^{zP_+}
\sinh zP_+\otimes
e^{zP_+} (D-C))=0.
\eea
Both equations can be summarized by the expression
\be
\exp(\pm 2z\sinh(z\Delta P_+) g)   =\exp(\pm 2z\Delta P_+)
\ee
leading to the same result for the
function $g$:
\be
g=\frac {\Delta P_+}{\sinh z\Delta P_+} .\label{bbc}
\ee

Finally, an explicit check shows that the $R$--matrix
\be
R=\exp\left\{ K\wedge \sinh zP_+ \, \frac {z\Delta P_+}
{\sinh z\Delta P_+}\right\}
\label{bbd}
\ee
also verifies property (\ref{aa}) for the last generator $K$.  Note
that the functions $f$ (\ref{bm}) and $g$ (\ref{bbc}) commute.

The result (\ref{bbd}) is in fact similar to a universal $R$--matrix
given in \cite{AA} for a Hopf algebra $\{ v,h\}$ which is isomorphic
to the Hopf subalgebra $\{K,P_+\}$ of $U_ziso(1,1)$. Therefore, since
the $R$--matrix (\ref{bbd}) does not depend on $P_-$ the universality
holds and it satisfies the quantum YBE (\ref{ac}). On the other hand,
it is worth remarking that expression (\ref{bbd})  has been also
obtained in \cite{Iran} following a different procedure.

\sect{Construction of $U_ziso(1,1)\oplus U_{-z}iso(1,1)$}

Let us consider two copies of $iso(1,1)$ with generators
$\{K^l,P_\pm^l\}$ $(l=1,2)$. The set of generators defined by
\be
J_3=K^1+K^2,\quad J_\pm= P_\pm^1 +P_\pm^2,\quad
N_3=K^1-K^2,\quad N_\pm= P_\pm^1 -P_\pm^2,\label{ca}
\ee
closes the algebra $t_4(so(1,1)\oplus so(1,1))\simeq iso(1,1) \oplus
iso(1,1)$. The  formal transformation (equivalent
 to a graded contraction \cite{Beyond})  defined by
\be
 (J,P_1,P_2,C_1,C_2,D)  :=
(\sqrt{\mu}N_3/2, J_+,
\sqrt{\mu}N_+,-J_-,\sqrt{\mu}N_-,J_3/2),
\label{cb}
\ee
gives rise to the following non-vanishing commutation relations
\bea
&& [J ,P_1]=  P_2, \qquad   [J,P_2]=  \mu  P_1, \qquad [D,P_i]=P_i, \cr
&& [J ,C_1]=  C_2, \qquad   [J,C_2]=  \mu  C_1, \qquad [D,C_i]=-C_i .
\label{cc}
\eea
For $\mu$ equal to $+1$, $0$ and $-1$ we obtain, in this order, the
commutators of $t_4(so(1,1)\oplus so(1,1))$,
$iiso(1,1)$ and $t_4(so(2)\oplus so(1,1))$. We denote these  three
algebras by $g_{\mu}$.

In the following, we will show how
the results presented in the previous section for the quantum
non-standard (1+1) Poincar\'e algebra provide a quantum structure for
the algebras $g_{\mu}$ and for the groups $G_{\mu}$ as well as their
universal $R$--matrices.

The invariance of $U_ziso(1,1)$ under the transformation  $z\to -z$
allows us to write $U_zt_4(so(1,1)\oplus so(1,1)) =U_ziso(1,1)\oplus
U_{-z}iso(1,1)$. The contraction (\ref{cc}) is implemented to the
quantum case by considering the following definition of the
contracted generators and   deformation parameter
\be
 (J,P_1,P_2,C_1,C_2,D;w)  :=
(\sqrt{\mu}N_3/2, J_+,
\sqrt{\mu}N_+,-J_-,\sqrt{\mu}N_-,J_3/2;z/\sqrt{\mu}),
\label{cbb}
\ee
where $w$ is the new (contracted) quantum parameter.
In this way we obtain the Hopf structure of $U_w g_\mu$. We omit the
explicit expressions so obtained since they are exactly the  quantum
algebras  $U_w^{(n)} g_{(\mu_1,0,+)}$
(with $\mu\equiv \mu_1$) given in \cite{Beyond}.

\subsect{Poisson--Hopf structure of $Fun(G_\mu)$}

Before getting the  quantum groups associated to $U_w g_\mu$ we
first study the algebra  $Fun(G_\mu)$ of smooth functions on
the group $G_\mu$.

A matrix realization of $g_\m$  in terms of $4\times 4$
real matrices is:
$$
J=\left(\begin{array}{cccc}
 0 & -\mu & 0 &0  \\ -1 & 0 &  0&
0 \\ 0 & 0 & 0 & 0\\ 0 & 0 & 0 & 0
\end{array}\right),\
 P_1 =\left(\begin{array}{cccc}
 0 & 0 & 0 &0  \\ 0 & 0 &  0&
 0 \\ 0 & 1 & 0 &0\\ 0 & 1 & 0 & 0
\end{array}\right), \
 P_2 =\left(\begin{array}{cccc}
 0 & 0  &0 &0 \\ 0 & 0 &  0&
0   \\1 & 0 & 0 & 0\\1 & 0 & 0 & 0
\end{array}\right),\nonumber
$$
\be
 D =\left(\begin{array}{cccc}
 0 &0 &0 &0 \\ 0 & 0 &  0&
0   \\ 0 & 0 & 0 & 1\\ 0 & 0 & 1 & 0
\end{array}\right),\
 C_1 =\left(\begin{array}{cccc}
 0 & 0 & 0 &0  \\ 0 & 0 &  0 & 0 \\ 0 &-1 & 0 & 0\\ 0 & 1 & 0 & 0
\end{array}\right), \
 C_2 =\left(\begin{array}{cccc}
 0 & 0 & 0 &0  \\ 0 & 0 &  0&
0   \\ -1 & 0 & 0 & 0\\ 1 & 0 & 0 & 0
\end{array}\right).\label{matrix}
\ee
Hence, a real  $4\times 4$ representation of the element
$g=e^{c_1C_1}e^{c_2C_2} e^{p_1P_1}e^{p_2P_2}e^{dD}e^{\theta J}\in G_\mu$
is given by \be
g=\left(\begin{array}{cccc}
 \Ck_{-\mu}(\theta) & -\mu \Sk_{-\mu}(\theta)& 0 &0  \\
 -\Sk_{-\mu}(\theta) &
\Ck_{-\mu}(\theta) &  0& 0 \\
t_{31} & t_{32} & \cosh d & \sinh d\\ t_{41} & t_{42} & \sinh d &
\cosh d \end{array}\right),\label{elem}
\ee
with
\bea
&&t_{31}=(p_2 -c_2) \Ck_{-\mu}(\theta) -(p_1 - c_1) \Sk_{-\mu}(\theta),\cr
&&t_{32}=(p_1 -c_1) \Ck_{-\mu}(\theta)  -\mu (p_2 - c_2)
\Sk_{-\mu}(\theta),\cr &&t_{41}=(p_2 +c_2) \Ck_{-\mu}(\theta) -(p_1 + c_1)
\Sk_{-\mu}(\theta),\\ &&t_{42}=(p_1 +c_1) \Ck_{-\mu}(\theta) -\mu (p_2 +
c_2)
 \Sk_{-\mu}(\theta) .\nonumber
\eea
The generalized sine and cosine functions are defined by
\be
\Ck_{{-\m}}(\theta)=\frac{e^{\sqrt{{\m}} \theta}+ e^{-\sqrt{{\m}}
\theta}}{2}, \qquad  \Sk_{{-\m}}(\theta)=\frac{e^{\sqrt{{\m}}
\theta}-e^{-\sqrt{{\m}} \theta}} {2\sqrt{{\m}}}.
\label{cd}
\ee
Note that for $\mu$ equal to $+1$ and $-1$ we recover  the hyperbolic
and elliptic trigonometric functions. The case $\mu=0$ corresponds to a
contraction of the group representation (\ref{elem}): $\Ck_0(\theta)=1$ and
$\Sk_0(\theta)=\theta$.

\begin{proposition}
The  fundamental  Poisson brackets
\bea
&&\{d,p_1\}= w\, \mu\, e^d
\Sk_{-\mu}(\theta),\qquad\qquad \{p_1,c_1\}=   w\,  \mu\, c_2,\cr
&&\{d,p_2\}= w\, ( e^d \Ck_{-\mu}(\theta)-1),
\qquad \{p_1,c_2\}=   w\,   c_1,\cr
&&\{\theta,p_1\}=w \,( e^d \Ck_{-\mu}(\theta)-1),
\qquad \{p_2,c_1\}= - w \,   c_1,
\label{ccc}\\
&&\{\theta,p_2\}= w\,   e^d \Sk_{-\mu}(\theta),
\qquad\qquad \{p_2,c_2\}= - w\,    c_2,
\nonumber
\eea
endow $Fun(G_\mu)$ with a Poisson--Hopf algebra structure.
\end{proposition}

The Poisson brackets (\ref{ccc}) comes from the Sklyanin bracket
induced from a
 classical $r$--matrix
\be
\{\Psi,\Phi\}= r^{\alpha\beta}\left(X_{\alpha}^L \Psi X_{\beta}^L
\Phi -  X_{\alpha}^R \Psi X_{\beta}^R \Phi \right)\label{ce}\quad
 \Psi,\Phi\in Fun(G_\mu).
\ee
In our case the $r$--matrix which satisfies the classical YBE  is given
by \be
 r=w\,(J\wedge P_1 +D\wedge P_2)\label{cf},
\ee
while left and right invariant vector fields are  deduced from
   (\ref{elem})
 \bea
&&
X_J^L=\partial_\theta,\qquad\qquad X_D^L=\partial_d,\cr
&&
X_{P_1}^L=e^d \Ck_{-\mu}(\theta)\partial_{p_1} +
e^d \Sk_{-\mu}(\theta)\partial_{p_2},\cr
&&
X_{P_2}^L=  e^d \Ck_{-\mu}(\theta)\partial_{p_2} +
\mu e^d \Sk_{-\mu}(\theta)\partial_{p_1},\label{cg}\\
&&
X_{C_1}^L=e^{-d} \Ck_{-\mu}(\theta)\partial_{c_1} +
e^{-d} \Sk_{-\mu}(\theta)\partial_{c_2},\cr
&&
X_{C_2}^L=e^{-d} \Ck_{-\mu}(\theta)\partial_{c_2} +
\mu e^{-d} \Sk_{-\mu}(\theta)\partial_{c_1},\nonumber
\eea
\bea
&&
X_J^R=\partial_\theta + \mu p_2 \partial_{p_1} + p_1 \partial_{p_2}
+\mu c_2 \partial_{c_1} + c_1 \partial_{c_2},\cr
&&
X_D^R=\partial_d +   p_1 \partial_{p_1} + p_2 \partial_{p_2}
-  c_1 \partial_{c_1} - c_2 \partial_{c_2},\label{ch}\\
&&
X_{P_1}^R=\partial_{p_1},\quad X_{P_2}^R=\partial_{p_2},\quad
X_{C_1}^R=\partial_{c_1},\quad X_{C_2}^R=\partial_{c_2}.\nonumber
\eea

\subsect{Hopf structure of $Fun_w(G_\mu)$}

Now we proceed to quantize the Poisson--Hopf algebra $Fun(G_\mu)$.
First we consider two sets of quantum coordinates
$\{\hat \te^l,\hat a_+^l,\hat a_-^l \}$ $(l=1,2)$ of
$Fun_z(ISO(1,1))$ for $l=1$ and of $Fun_{-z}(ISO(1,1))$ for $l=2$.
Then we construct the Hopf algebra
$Fun_z(ISO(1,1))\oplus Fun_{-z}(ISO(1,1))$ with the results of  Prop.\
2 and by using the new coordinates defined by
\be
\hat a=\frac 12(\hat \te^1 +\hat \te^2),\quad   \hat a_\pm=\frac
12(\hat a_\pm^1 +\hat a_\pm^2),\quad\hat  b=\frac 12(\hat \te^1 -\hat
\te^2),\quad \hat  b_\pm=\frac 12(\hat a_\pm^1 -\hat
a_\pm^2).\label{ci}  \ee
Next we apply the quantum contraction induced at the group level from
(\ref{cbb}):
\be
 (\hat \theta,\hat p_1,\hat p_2,\hat c_1,\hat c_2,\hat d;w)  :=
( 2\hat  b/\sqrt{\mu},\hat  a_+,
\hat b_+/\sqrt{\mu},-\hat a_-,\hat b_-/\sqrt{\mu},2\hat
a;z/\sqrt{\mu}), \label{cj}
\ee
obtaining in this way the quantization of $Fun(G_\mu)$.
The final result is summarized as follows:

\begin{proposition}
 The  Hopf algebra
$Fun_w(G_\mu)$  is given  by non-vanishing commutators
\bea
&&[\hat d,\hat p_1]= w\, \mu\, e^{\hat d}  \Sk_{-\mu}(\hat
\theta),\qquad\qquad [\hat p_1,\hat c_1]=   w\,\mu\,\hat  c_2,\cr
&&[\hat d,\hat p_2]= w\, ( e^{\hat d} \Ck_{-\mu}(\hat \theta)-1),\qquad
[\hat p_1,\hat c_2]=   w\,\hat  c_1,\cr
&&[\hat \theta,\hat p_1]=w\, ( e^{\hat d} \Ck_{-\mu}(\hat \theta)-1),\qquad
[\hat p_2,\hat c_1]= - w \,\hat  c_1, \label{ck}\\
&&[\hat \theta,\hat p_2]= w\,   e^{\hat d} \Sk_{-\mu}(\hat \theta),
\qquad\qquad [\hat p_2,\hat c_2]= - w \,\hat  c_2; \nonumber
\eea
coproduct, counit and antipode
\bea
&&\Delta(\hat \theta)=\hat \theta\otimes 1 + 1 \otimes\hat  \theta,
\qquad \Delta(\hat d)=\hat d\otimes 1 + 1 \otimes\hat  d,\cr
&&
\Delta(\hat p_1)=\hat p_1\otimes 1 + e^{\hat d}\Ck_{-\mu}(\hat \theta)
\otimes\hat  p_1
+ \mu e^{\hat d} \Sk_{-\mu}(\hat \theta)\otimes\hat  p_2,\cr
&&
\Delta(\hat p_2)=\hat p_2\otimes 1 + e^{\hat d}
\Ck_{-\mu}(\hat \theta)\otimes\hat  p_2 +
e^{\hat d} \Sk_{-\mu}(\hat \theta)\otimes\hat  p_1,\label{cl}\\
&&
\Delta(\hat c_1)=\hat c_1\otimes 1 + e^{-\hat d}
\Ck_{-\mu}(\hat \theta)\otimes\hat  c_1 + \mu
e^{-\hat d} \Sk_{-\mu}(\hat \theta)\otimes\hat  c_2,\cr
&&
\Delta(\hat c_2)=\hat c_2\otimes 1 + e^{-\hat d}
\Ck_{-\mu}(\hat \theta)\otimes\hat  c_2 +
e^{-\hat d} \Sk_{-\mu}(\hat \theta)\otimes\hat  c_1;\nonumber
\eea
\be
\epsilon(X) =0,\qquad
X\in\{\hat \theta,\hat p_i,\hat c_i,\hat d\};
\label{cch}
\ee
\bea
&& \gamma(\hat \theta)=  -\hat  \theta,\qquad
\gamma(\hat  d)= -\hat  d,\cr
&&
 \gamma(\hat  p_1)=  - e^{-\hat d}\Ck_{-\mu}(\hat \theta)\hat p_1 - \mu
e^{-\hat d} \Sk_{-\mu}(\hat \theta)\hat  p_2,\cr
&&
\gamma(\hat  p_2)= - e^{-\hat d} \Ck_{-\mu}(\hat \theta) \hat  p_2 -
e^{-\hat d} \Sk_{-\mu}(\hat \theta) \hat  p_1,\label{cm}\\
&&
 \gamma(\hat  c_1)=- e^{\hat d} \Ck_{-\mu}(\hat \theta)\hat  c_1 - \mu
e^{\hat d} \Sk_{-\mu}(\hat \theta)\hat  c_2,\cr
&&
\gamma(\hat  c_2)= - e^{\hat d} \Ck_{-\mu}(\hat \theta)\hat  c_2 -
e^{\hat d} \Sk_{-\mu}(\hat \theta)\hat   c_1.\nonumber
\eea
\end{proposition}

The final step in this quantization process consists on deducing  the
universal $R$--matrix for $U_w g_\mu$.
We write two $R$--matrices (\ref{bbd}) $R_{ z}^1$ and $R_{- z}^2$ with
generators $\{K^l,P^l_\pm\}$ $(l=1,2)$ for   the two copies
$U_{\pm z} iso(1,1)$ and
compute the product ${\cal R}=R^1_z R^2_{-z}$:
\bea
&&
\back {\cal R}=
\exp\left\{ K^1\wedge \sinh zP^1_+ \,\frac {z\Delta P^1_+}
{\sinh z\Delta P^1_+}\right\}
\exp\left\{- K^2\wedge \sinh zP^2_+ \, \frac {z\Delta P^2_+} {\sinh
z\Delta P^2_+}\right\} \cr
 &&
\ = \exp\left\{ (K^1\wedge \sinh zP^1_+ \, \Delta P^1_+
\sinh z\Delta P^2_+\right.\cr
&&
\qquad\quad \left. - K^2\wedge \sinh zP^2_+ \,
 \Delta P^2_+\sinh z\Delta P^1_+)
\, \frac{z}{\sinh z\Delta P^1_+\sinh z\Delta P^2_+}\right\} .
\eea
We introduce the change of generators
(\ref{ca}) and, afterwards, we apply the quantum contraction
(\ref{cbb}). The  final expression for the universal $R$--matrix of
$U_w g_\mu$ (denoted by ${\cal R}_w$) is
\be
{\cal R}_w=\exp\{(M_1 N_1 + M_2 N_2) L \},\label{cn}
\ee
where
\bea
&&\back\back
M_1=D\wedge \Ck_{-\mu}(wP_1/2)\sinh(wP_2/2) +  J\wedge
\Sk_{-\mu}(wP_1/2) \cosh(wP_2/2),\cr
&&\back\back
M_2=\mu D\wedge \Sk_{-\mu}(wP_1/2)\cosh(wP_2/2) +  J\wedge
\Ck_{-\mu}(wP_1/2) \sinh(wP_2/2),\cr
&&\back\back
N_1=\mu \Delta P_1 \Sk_{-\mu}(w\Delta P_1/2) \cosh(w\Delta P_2/2) -
\Delta P_2 \Ck_{-\mu}(w\Delta P_1/2) \sinh(w\Delta P_2/2),\cr
&&\back\back
N_2= \Delta P_2 \Sk_{-\mu}(w\Delta P_1/2) \cosh(w\Delta P_2/2) -
\Delta P_1 \Ck_{-\mu}(w\Delta P_1/2) \sinh(w\Delta P_2/2),\cr
&&\back\back
L=\frac{2w}{\Ck_{-\mu}(w\Delta P_1) -\cosh(w\Delta P_2)}.\label{co}
\eea

An interesting  idea  naturally arising from this result
would be the use of the
FRT construction \cite{FRT} to quantize $Fun(G_\mu)$. In fact,
the matrix representation
(\ref{matrix}) substituted in (\ref{cn}) gives rise to a
particular representation of
${\cal R}_w$:
\bea
&&{\cal R}_w=\exp\{w r\} = \exp\{w(J\wedge P_1 + D\wedge P_2)\}\cr
&&\qquad = I\otimes I + w(J\wedge P_1 + D\wedge P_2)  +\mu w^2
P_1\otimes P_1 , \eea
where $I$ is the four dimensional  identity matrix.  In this
representation the commutation rules of the group coordinates
$(\hat d,\hat \theta,\hat p_i,\hat c_i)$ would be deduced from
equation
\be
{\cal R}_w T_1 T_2 = T_2 T_1 {\cal R}_w ,
\ee
  $T$ is the generic element of the group $G_\mu$ (\ref{elem}),
 $T_1=T\otimes
I$ and $T_2=I\otimes T$. Lengthy computations show that
commutators so obtained are
exactly those given in  (\ref{ck}) up to a global change of sign
 in the deformation
parameter $w$. Furthermore, coproduct  (\ref{cl}), counit (\ref{cch})
 and antipode
(\ref{cm}) can be got from relations $\Delta(T)=T\dot\otimes T$,
$\epsilon(T)=I$ and $\gamma(T)=T^{-1}$.

\sect{Universal quantizations of Weyl subalgebras}

The Lie brackets
\be
[J,P_1]=P_2,\qquad [J,P_2]=\mu P_1,\qquad [P_1,P_2]=0,\label{da}
\ee
correspond   for $\mu$ negative, positive and zero to the two dimensional
Euclidean,  Poincar\'e and Galilei algebras, respectively.
We can enlarge these algebras by means of a dilation generator $D$:
\be
[D,P_i]=P_i,\qquad [D,J]=0.\label{db}
\ee
These enlarged algebras are the algebras of similarities of the
Euclidean, Minkowskian or Galilean planes, and will be denoted by
$s_\mu$. They are naturally the ``Weyl" subalgebras of the corresponding
conformal algebras in two dimensions.
Although the  conformal algebras of the  family
(\ref{da}) are  indeed $so(3,1)$,
$iso(2,1)$ and $so(2,2)$ for $\mu<,=,>0$ \cite{Beyond},    the crucial
point  is that each of the algebras in the family $g_\mu$ also contains a
subalgebra isomorphic to the Weyl subalgebra of the full conformal algebras
of the Euclidean,  Poincar\'e and Galilei spaces. Moreover, the Hopf algebra
$U_wg_\mu$   preserves this property, that is,   $U_wg_\mu$ includes quantum
Weyl subalgebras that deform (\ref{da},\ref{db}).

\begin{proposition}
The algebras $U_w s_\mu$
given by \cite{Beyond}:
\bea
&&\back {\Delta P_1 =1 \otimes P_1 + P_1\otimes 1,
   \qquad \Delta P_2 =1 \otimes P_2 + P_2\otimes 1,} \cr
&& \back {\Delta J =
   e^{-\tfrac w2 P_2 }\Ck_{-\m_1}(wP_1/2) \otimes J  +
   J \otimes \Ck_{-\m_1}(wP_1/2) e^{\tfrac w2 P_2 }} \cr
&& \quad
  -e^{-\tfrac w2 P_2}\Sk_{-\m_1}(wP_1/2) \otimes\m_1 D
  + \m_1 D \otimes \Sk_{-\m_1}(wP_1/2) e^{\tfrac w2 P_2 },\cr
&&\back  {\Delta D =
   e^{-\tfrac w2 P_2 }\Ck_{-\m_1}(wP_1/2) \otimes D  +
   D \otimes \Ck_{-\m_1}(wP_1/2) e^{\tfrac w2 P_2 }} \cr
&& \quad
  -e^{-\tfrac w2 P_2}\Sk_{-\m_1}(wP_1/2) \otimes J
  + J \otimes \Sk_{-\m_1}(wP_1/2) e^{\tfrac w2 P_2 };\nonumber
\eea
\be
\epsilon(X) =0; \qquad
  \gamma(X) = -e^{{w}  P_2}\ X\ e^{-{w}  P_2},
  \qquad X \in \{ J,P_i,D\} ;
\ee
\bea
&&[J,P_1]=\frac 2w\sinh(wP_2/2)\Ck_{- \m_1}(wP_1/2),\cr
&&[J,P_2]=\frac 2w\m_1\Sk_{- \m_1}(wP_1/2)\cosh(wP_2/2),\cr
&&[D,P_1]=\frac 2w \Sk_{- \m_1}(wP_1/2)\cosh(wP_2/2),\cr
&&[D,P_2]=\frac 2w\sinh(wP_2/2)\Ck_{- \m_1}(wP_1/2),\cr
&&[P_1,P_2]=0,\qquad [D,J]=0;
 \nonumber
\eea
  are quasitriangular Hopf algebras with universal $R$--matrix
(\ref{cn},\ref{co}).
\end{proposition}

Furthermore, it is   clear that by taking the generators
$C_i\equiv 0$ and the group parameters $\hat c_i\equiv 0$ in Props.\ 3
and 4 we find a Poisson--Hopf algebra structure for $Fun(S_\mu)$ and a
quantum Hopf algebra $Fun_w(S_\mu)$.

\sect{Concluding remarks}

A combined approach of
the construction $U_z A\oplus U_{-z} A$ ($A$ being either  a Lie
algebra or the algebra of functions on the Lie group) together with a
quantum contraction   provide  a simultaneous universal quantization
for the algebras $g_\mu$ in the family (\ref{cc}).

One of the groups in the family $G_\mu$  can be realized as a kinematical
group: the group $G_{-1}\equiv T_4(SO(2)\otimes SO(1,1))$ is
isomorphic to the (2+1) expanding Newton--Hooke group \cite{LLB}, the motion
group of a   non-relativistic space-time with    constant   negative
curvature. Time is absolute in such a universe   and a space-time
contraction leads from it to the Galilean case. An adapted basis for $G_{-1}$
is formed by a time translation $\tilde H$, two spatial translations $\tilde
P_i$, two boosts $\tilde K_i$ and one spatial rotation $\tilde J$, with
corresponding group coordinates $\{t,x_i,v_i,\psi\}$ $(i=1,2)$.
 All expressions
obtained for $G_{-1}$ in   section 3 can be written in terms of  these
new generators and group coordinates by means of the following
isomorphisms:  \bea
&&
\tilde J\equiv J,\quad \tilde P_i\equiv\frac 12(P_i+C_i),\quad
\tilde K_i\equiv\frac 12(P_i-C_i),\quad \tilde H\equiv -D,\\
&&
  \psi \equiv \theta,\quad   x_i\equiv 2(p_i+c_i),\quad
  v_i\equiv 2(p_i-c_i),\quad  t\equiv -d.
\eea

An open problem to be solved is the construction of a universal
$R$--matrix for the non-standard quantum deformation of $sl(2,\R)$
which would provide  a set of
universal $R$--matrices for the whole set of (2+1) non-standard quantum
algebras, following the method just described. We recall that among them
there are some rather interesting cases from a physical point of view: the
conformal algebras of the (1+1) Poincar\'e ($so(2,2)$) and two
dimensional Euclidean spaces ($so(3,1)$), besides a ``null-plane" (2+1)
Poincar\'e algebra.

\bigskip
\bigskip

\noindent
{\large{{\bf Acknowledgements}}}

\bigskip

This work has been partially supported by a DGICYT project
(PB92--0255) from the Ministerio de Educaci\'on y Ciencia de
Espa\~na and by an Acci\'on Integrada Hispano--Italiana (HI--059).

\end{document}